# Numerical study on transient heat transfer under soil with plastic mulch in agriculture applications using a nonlinear finite element model

Carlos Armando De Castro and Orlando Porras Rey

**Abstract**

In this paper is developed a simple mathematical model of transient heat transfer under soil with plastic mulch in order to determine with numerical studies the influence of different plastic mulches on the soil temperature and the evolutions of temperatures at different depths with time. The governing differential equations are solved by a Galerkin Finite Element Model, taking into account the nonlinearities due to radiative heat exchange between the soil surface, the plastic mulch and the atmosphere. The model was validated experimentally giving good approximation of the model to the measured data. Simulations were run with the validated model in order to determine the optimal combination of mulch optical properties to maximize the soil temperature with a Taguchi's analysis, proving that the material most used nowadays in Colombia is not the optimal and giving quantitative results of the properties the optimal mulch must possess.

**Keywords:** heat equation, finite elements, nonlinear boundary conditions, plastic mulch, atmospheric radiation, meteorological inputs

## SYMBOLS

- $\tau_l$     Mulch transmittance to long wave radiation
- $\rho_l$     Mulch reflectivity to long wave radiation
- $\tau_s$     Mulch transmittance to solar radiation
- $\rho_s$     Mulch reflectivity to solar radiation
- $\epsilon_m$     Mulch emissivity
- $\epsilon_s$     Soil surface emissivity
- $a_s$     Soil albedo
- $k$     Soil thermal conductivity [W/m K]
- $\rho c_p$     Soil specific heat per volume unit [J/m³ K]
- $h_i$     Convective heat transfer coefficient between the soil surface and the mulch [W/m² K]
- $h_o$     Convective heat transfer coefficient between the mulch and the ambient air [W/m² K]
- $R_s$     Direct and diffuse solar shortwave radiation [W/m²]
- $T_a$     Ambient air temperature [K]
- $v$     Wind velocity [m/s]
- $T_m$     Mulch temperature [K]
- $T$     Soil temperature [K]
- $q$     Heat flux [W/m²]
- $\sigma$     Boltzmann constant, $5.67 \times 10^{-8}$ W/m² K⁴

## 1 INTRODUCTION

Plastic mulch has had a great rise among plant growth techniques in the last decades. It consists in covering the soil surface with a plastic film that retains humidity and heat, so soil temperature is elevated and soil properties are improved so a better quality of grown plants is obtained and some unwanted parasites are killed by the high temperatures achieved [1]. In Colombia, plastic mulches are used in the flower's industry among others, but there are not investigations in the country dedicated to develop a general model to compute the temperature distribution in soil in order to determine the plastic mulch's properties needed to specific applications, so it is necessary to develop a model that helps make that choosing without trial-and-error procedures, and that helps visualize the influence of plastic mulch in the soil's thermal response under specific meteorological conditions, with a particular case being the protection of seeds in terrains in the Savanna of Bogotá where there are reached temperatures below zero at some times of the year.

## 2 GOVERNING EQUATIONS

It is considered that the heat absorbed by the soil surface is transmitted entirely to the bottom of the soil [2], [3], so the temperature distribution under soil is governed by the one-dimensional transient heat diffusion equation:

$$\rho c_p \frac{\partial T}{\partial t} - \frac{\partial}{\partial z}\left(k \frac{\partial T}{\partial z}\right) = 0 \qquad (1)$$

The boundary conditions of equation (1) are the net heat flux absorbed at the soil surface and the temperature at the bottom of the soil. The net heat flux absorbed at the soil surface and transmitted by conduction to the lower layers of the soil is the sum of the convection with the plastic film, the net solar radiation and the net infrared radiation:

$$-k \frac{\partial T(0)}{\partial z} = h_i(T_m - T) + \frac{\tau_s(1-a_s)}{1-\rho_s a_s} R_s + \frac{\epsilon_s \sigma}{1-\rho_l+\rho_l \epsilon_s}\left[\tau_l \epsilon_{sky} T_a^4 + \epsilon_m T_m^4 - (1-\rho_l)T^4\right] \qquad (2)$$

Sky emissivity is calculated by the Swinbank formula [4]:

$$\epsilon_{sky} = 9.2 \times 10^{-6} T_a^2 \qquad (3)$$

And the internal heat transfer coefficient was found experimentally by Garzoli and Blackwell [5] to be $h_i = 7.2 \text{ Wm}^{-2}\text{K}^{-1}$. The terms in equation (2) regarding radiation are found using a ray-diagram and geometric series [2]. The bottom of the soil is the depth at which there are no variations of temperature with time, that being a realistic assumption because soil is modeled as a semi-infinite body:

$$T(t, z_\infty) = T_\infty \qquad (4)$$

Equation (2) shows that the temperature of the mulch must be known at each time step, which is achieved by the heat balance on the mulch. Because of the small thickness of the plastic films generally used, it is considered that the mulch does not store energy, so the heat balance is:

$$h_i(T - T_m) + h_o(T_a - T_m) + \left[1 - \rho_s - \frac{\tau_s(1 - a_s + \tau_s a_s)}{1 - \rho_s a_s}\right] R_s + \left[1 - \rho_l - \frac{\tau_l(\tau_l + \epsilon_s(1 - \tau_l))}{1 - \rho_l + \rho_l \epsilon_s}\right] \epsilon_{sky} \sigma T_a^4 -$$
$$\left[2 - \frac{(1 - \epsilon_s)(1 - \tau_l - \rho_l)}{1 - \rho_l + \rho_l \epsilon_s}\right] \epsilon_m \sigma T_m^4 + \left[1 - \frac{\tau_l + \rho_l}{1 - \rho_l + \rho_l \epsilon_s}\right] \epsilon_s \sigma T^4 = 0 \qquad (5)$$

The external convective heat transfer coefficient was found experimentally by Garzoli and Blackwell [5] to be:

$$h_o = 7.2 + 3.8v \text{ Wm}^{-2}\text{K}^{-1} \qquad (6)$$

Like in equation (2), the terms in equation (5) regarding radiation are found using a ray-diagram and geometric series [2].

## 3   FINITE ELEMENT MODEL

Using a uniform mesh of elements each one of length $h_e$ and quadratic Lagrange interpolation functions, the weak form of equation (1) in explicit matrix form using the Galerkin method and with element-wise constant values of the soil thermal properties is [8]:

$$\frac{\rho c_p}{30} \begin{bmatrix} 4 & 2 & -1 \\ 2 & 16 & 2 \\ -1 & 2 & 4 \end{bmatrix} \begin{Bmatrix} \dot{T}_1^e \\ \dot{T}_2^e \\ \dot{T}_3^e \end{Bmatrix} + \frac{k}{3h_e} \begin{bmatrix} 7 & -8 & 1 \\ -8 & 16 & -8 \\ 1 & -8 & 7 \end{bmatrix} \begin{Bmatrix} T_1^e \\ T_2^e \\ T_3^e \end{Bmatrix} = \begin{Bmatrix} q_1 \\ 0 \\ -q_3 \end{Bmatrix} \qquad (7)$$

After assembly of the element matrices results the global matrix form of the finite element model:

$$[C]\{\dot{T}\} + [K]\{T\} = \{F\} \qquad (8)$$

Approximating the time derivative by a backward finite-differences model, we have:

$$[\hat{K}]\{T\}^{n+1} = \{\hat{F}\}^{n+1} \qquad (9)$$

Where $[\hat{K}] = [C] + \Delta t[K]$ and $\{\hat{F}\}^{n+1} = [C]\{T\}^n + \Delta t\{F\}^{n+1}$, and the superscript represents the time step. Equation (10) is a nonlinear algebraic system of equations because of the 4th power of the temperature terms in the boundary conditions given by equation (2). The system is solved by the fixed point iteration [8]:

$$\{T\}_{p+1}^{n+1} = [\widehat{K}]^{-1}\{\widehat{F}\}_p^{n+1} \tag{10}$$

Additionally, equation (5) is also solved simultaneously by fixed point iteration in order to know the mulch temperature and have a complete solution of all of the variables involved. The solution was implemented in an algorithm in MATLAB.

## 4  MODEL VALIDATION

A mulch of low density polyethylene (LDPE) was constructed in a terrain at the Savanna of Bogotá. The mulch optical properties and the soil thermal properties were characterized and then meteorological data and soil surface temperature were measured simultaneously using the portable automated weather station Casella NOMAD and an OMEGA 871A NiCr-NiAl thermocouple. The spectrophotometer NICOLET 380 FTIR was used to characterize the used LDPE properties to long wave radiation and the spectrophotometer Perkin Elmer Lambda 3 UV/VIS was used to characterize the used LDPE mulch properties to solar shortwave radiation. According to Wien's Displacement Law, the infrared wavelengths of interest in this case are in the range of 8700 nm to 10600 nm, also, the solar radiation sensor of the weather station measures in the range of 400 nm to 1100 nm. The optical properties of the used LDPE in those ranges were found by numerical integration and are shown in Table 1.

**Table 1.** Measured optical properties of the LDPE mulch.

| Property | Value |
|---|---|
| $\tau_l$ | 0.6 |
| $\rho_l$ | 0.398 |
| $\tau_s$ | 0.733 |
| $\rho_s$ | 0.265 |

The soil thermal properties were calculated using a 75 W electric resistance and measuring the variation of temperature in a small sample of soil in a range of 10 min. The estimated soil thermal properties are thermal conductivity of 2.2 W/m×K and specific heat per unit volume of $1.01 \times 10^6$ J/m$^3$.

Data were measured each 5 min between 09:45 and 14:40 of December 5$^{th}$ of 2008. The meteorological input data measured by the automated weather station is shown in Figures 1 to 3.

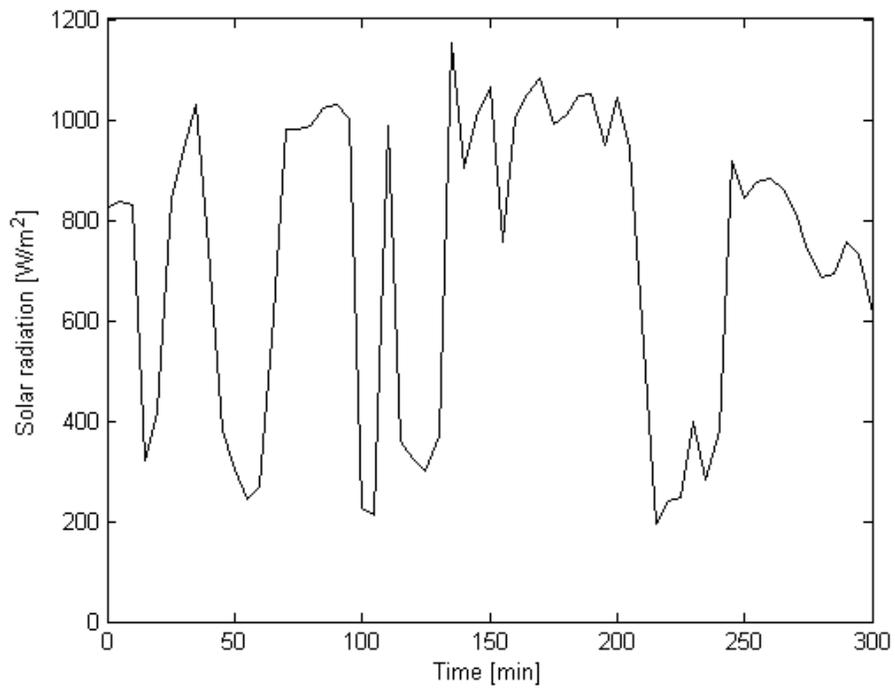

**Figure 1.** Solar radiation.

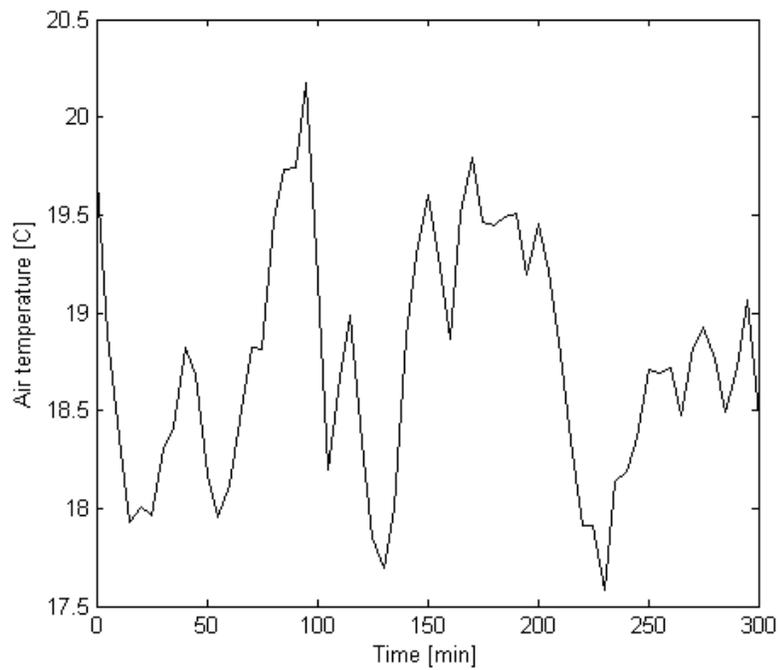

**Figure 2.** Ambient air temperature.

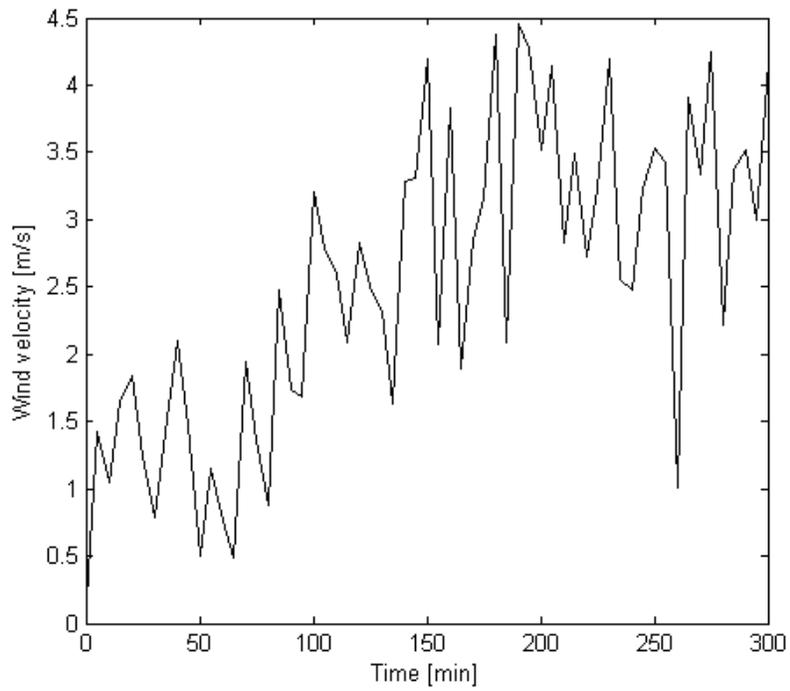

**Figure 3.** Wind velocity.

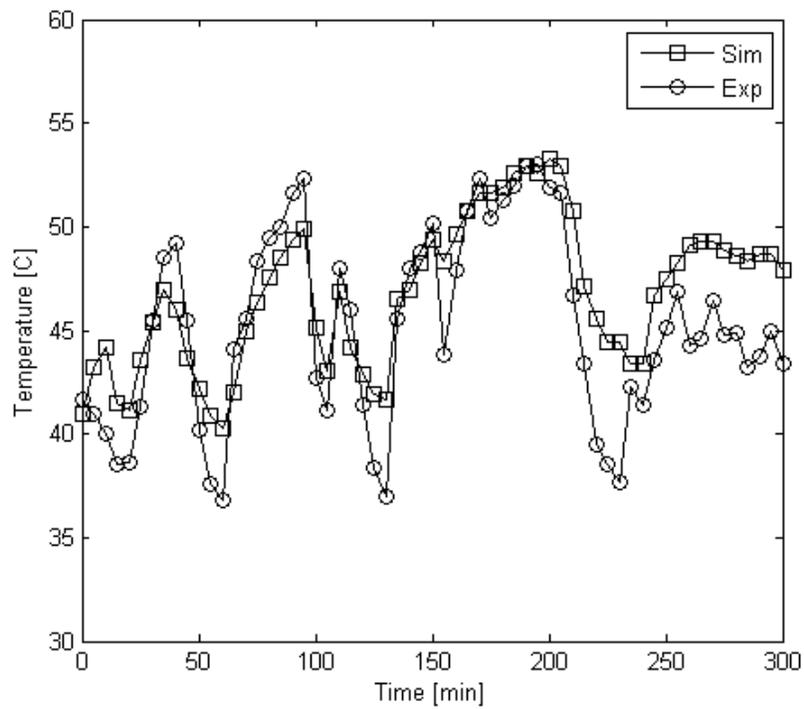

**Figure 4.** Comparison between measured (Exp) and simulated (Sim) soil surface temperature.

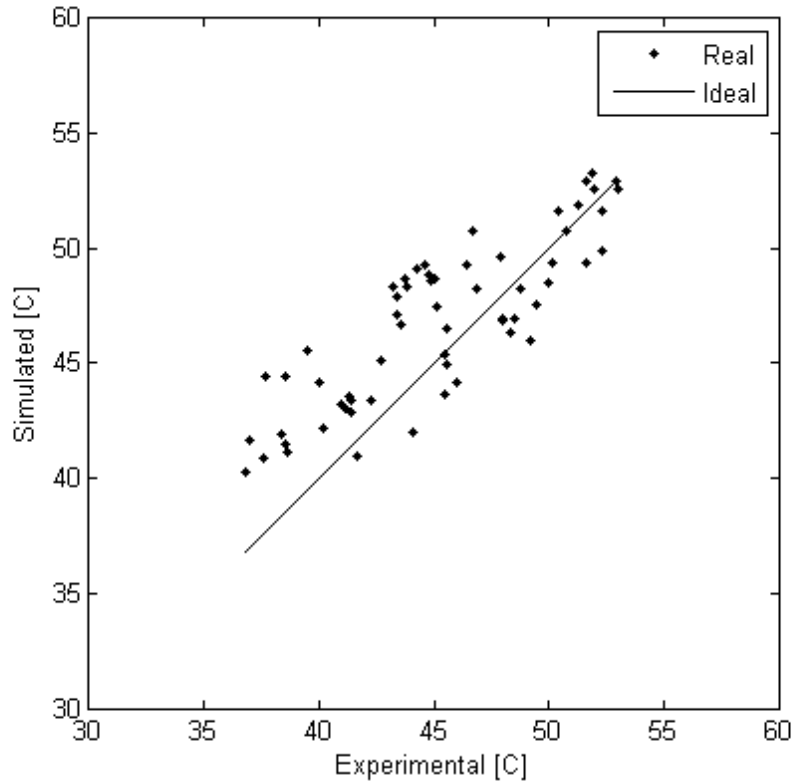

**Figure 5.** Correlation between experimental and simulated temperatures.

Simulated temperature was calculated using 50 elements and with a value for the bottom of the soil of 1 m. It can be seen in Figures 4 and 5 that there exists good agreement of the model and the measured values. It is worth noting that the weather station gives a 5 min average of meteorological data, so it is possible that the errors seen when large changes of solar radiation occur are due to the way data is acquired and not because of the model, that is better seen in Figure 4 after 200 min of experimentation, where rapid instantaneous changes of solar radiation were seen due to cloud movement. The mean relative error between measured and simulated soil surface temperature was 5.76%.

## 5    OPTIMIZATION OF MULCH'S PROPERTIES

With the validated model it was performed a Taguchi analysis [10] to find the optimal combination of properties to achieve the highest possible soil's temperature. The model was run with the input variables shown in Figures 1 to 3 and with the values of the soil's properties of the experimental validation. The combinations tested and the results of the maximum temperature are shown in Table 2. The Taguchi analysis and the optimal combination of properties found are shown in Figure 6 and Table 3.

**Table 2.** Simulations run for the Taguchi analysis.

| No. Sim. | PLASTIC PROPERTIES | | | | Tmax [°C] |
|---|---|---|---|---|---|
| | $\tau_s$ | $\rho_s$ | $\tau_l$ | $\rho_l$ | |
| 1 | 0.01 | 0.01 | 0.01 | 0.01 | 41.785 |
| 2 | 0.01 | 0.2 | 0.2 | 0.2 | 41.023 |
| 3 | 0.2 | 0.01 | 0.2 | 0.7 | 46.954 |
| 4 | 0.2 | 0.2 | 0.7 | 0.01 | 41.276 |
| 5 | 0.2 | 0.7 | 0.01 | 0.2 | 40.943 |
| 6 | 0.7 | 0.01 | 0.7 | 0.2 | 53.480 |
| 7 | 0.7 | 0.2 | 0.01 | 0.7 | 59.881 |

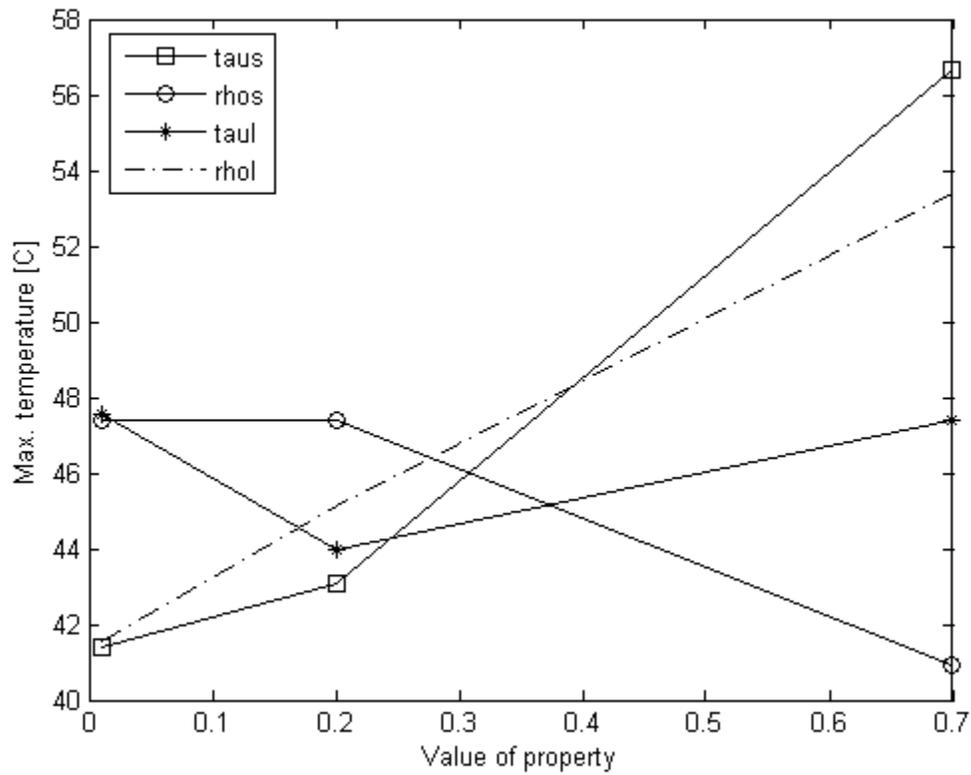

**Figure 6.** Taguchi analysis of mulch optical properties' influence in soil temperature.

**Table 3.** Optimal combination of mulch's properties to achieve maximum temperature.

| | |
|---|---|
| $\tau_s$ | >0.7 |
| $\rho_s$ | <0.2 |
| $\tau_l$ | <0.01 |
| $\rho_l$ | >0.7 |

## 6  SOIL'S TEMPERATURE VARIATION WITH DEPTH AND TIME

The model was run for a typical meteorological year of Bogotá, Colombia, to observe the variations of temperature with time and depth. It was observed that the daily fluctuations practically cease at a depth of 1.00 m, validating the first supposition made to compare the model with the validated data. However, the fluctuations for the whole year cease at a depth of 2.5 m. Figure 7 shows the temperatures at different depths compared with the ambient temperature for a random week of the simulated year.

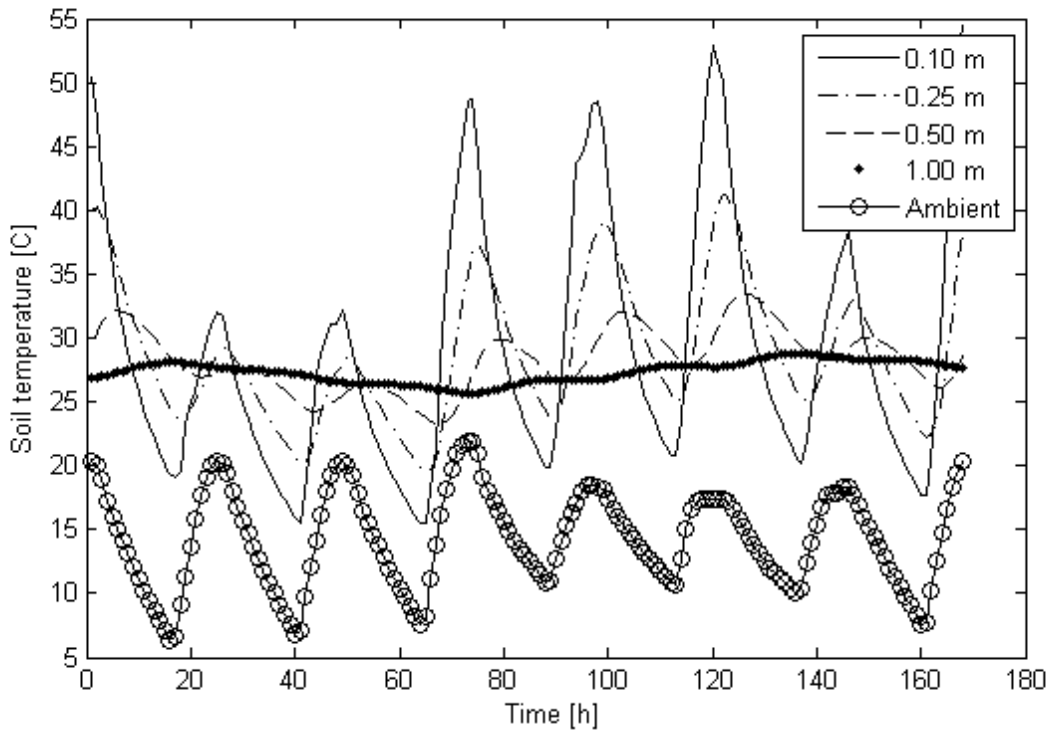

**Figure 7.** Simulated temperatures at different depths in comparison to ambient temperature for a week of simulation.

## 7  CONCLUSIONS

The one-dimensional nonlinear finite element model developed has good agreement with experimental data, so it is inferred that it is useful to predictions of temperature of soil with different kinds of plastic mulches, and will help to develop curves of mulch selection to specific applications. Results from the Taguchi analysis shows that the material used for plastic mulch in Colombia (LDPE) is not the optimal. The future work will be pointed to develop non-dimensional curves to help select mulches to specific applications; another work will be done to analyze possible solutions regarding plastic mulches to protect the flower seeds when there are freezing temperatures at the Savanna of Bogotá.